\documentclass[a4paper]{article}
\usepackage[latin1]{inputenc}
\usepackage{graphics}

\usepackage{amssymb}
\usepackage{amsmath}
\usepackage[thmmarks, amsmath]{ntheorem}
\usepackage[all]{xy}

\def\A{{\mathcal A}}
\def\B{{\mathcal B}}

\def\H{{\mathcal H}}

\def\RR{{\mathbb R}}
\def\CC{{\mathbb C}}
\def\NN{{\mathbb N}}

\def\PP{{\mathbb P}}

\def\End{\mbox{\rm End}}

\def\bra{\langle}
\def\ket{\rangle}

\def\bea{\begin{eqnarray}}
\def\eea{\end{eqnarray}}

\def\be{\begin{equation}}
\def\ee{\end{equation}}

\newenvironment{rem}[1][{}]{\smallbreak \noindent  {\bf Remark #1}\small }

\newtheorem{theorem}{Theorem}
\newtheorem{definition}{Definition}
\newtheorem{lemma}{Lemma}
\newtheorem{cor}{Corollary}
\newtheorem{propo}{Proposition}
\theoremstyle{nonumberplain}
\theorembodyfont{\normalfont}
\theoremseparator{:}
\theoremsymbol{$\P$}

\newtheorem{demo}{Proof}
\begin{document}
\title{Estimating noncommutative distances on graphs}
\author{Fabien Besnard}

\maketitle
\begin{abstract}
We report on some findings concerning Connes' noncommutative distance $d$ on a weighted undirected graph $G$. Our main result is the lower bound $\ell/\Delta(G)\le d$ where $\ell$ is the geodesic distance and $\Delta(G)$ the   degree of $G$. It is obtained thanks to an auxiliary spectral triple on the collection of the edges of $G$.
\end{abstract}
\section{Introduction}
The fundamental intuition underlying Connes' noncommutative geometry is the identification of the inverse line element $ds^{-1}$ with the Dirac operator \cite{redbook}. On a Riemannian spin manifold $(M,g)$, it is possible to make this intuition precise through the formula
\begin{equation}
d_g(p,q)=\sup\{|f(p)-f(q)||\  \|[D,f]\|\le 1\}\label{formuledistance1}
\end{equation}
where $p,q\in M$ and $f$ is a smooth function. Here $d_g$ is the geodesic distance associated to the metric $g$, and $D$ is the canonical Dirac operator associated to the metric and spin structure. Formula \eqref{formuledistance1} can be generalized to a spectral triple $(\A,\H,D)$, giving rise to the definition of the noncommutative distance
\begin{equation}
d^D(\varphi,\psi)=\{|\varphi(a)-\psi(a)||\ a\in\A, \|[D,a]\|\le 1\}\label{formuledistance2}
\end{equation}
between two pure states $\varphi,\psi$ on the $C^*$-algebra $\A$. This generalized distance (it can be infinite), has been studied and computed in various contexts \cite{bimonte1994distances,atzmon1996distances,ikm,cagnache2010spectral,wallet2012connes}. It has also been extended to Lorentzian manifolds \cite{parfionov2000connes,francoprocfir,canarutto2019distance} and the quantization has been studied in a few simple cases \cite{rovnice,besnice}. We refer to \cite{martinetti2016monge} for a more thourough review. In this paper we will be exclusively concerned with the case where the algebra is commutative and finite. Then \eqref{formuledistance2} becomes a distance on a finite  graph $G$ with weights   given by the inverse   of the non-zero entries of  $D$. Though this might seem to be the simplest of all cases, it raises interesting questions.  The first of these questions is the relationship between $d^D$ and the geodesic distance $\ell^D$ defined by the weights. It is known that $\ell^D\not=d^D$ as soon as there are $3$ connected points or more in $G$, and that the inequality 
\begin{equation}
d^D\le \ell^D
\end{equation}
always hold. This is already a major difference with the smooth case, and it has led several authors to attach a different Dirac operator to a finite graph in order to recover the equality with the geodesic distance \cite{dimakis1998connes,requardt2002distance}. We will not follow this path in this paper, except at some intermediate stage, since our goal is to study $d^D$ for itself. There might several reasons to do so. First, when one approximates a smooth manifold using a lattice, the natural discretization of the Dirac operator is of the kind studied here \cite{bimonte1994distances,MVS}. At a more fundamental level, Dirac operators of this kind also arose in relation to quantum gravity \cite{rov1}.  Let us also mention that exploring the properties of $d^D$ is an interesting challenge in its own sake. For example, we have not been able to prove (or disprove) that $d^D$ grows as soon as one removes an edge, and this is only one of many gaps in our current understanding.

Explicit computations being out of reach as soon as there are as few as $4$ connected points in  a graph with generic weights, we will rather look for tools to estimate the noncommutative distance. We will give two such tools. The first one is to decompose $G$ into a sequence of ``blobs'' and chains and then estimate the distance through the lengths of the chains and the diameters of the blobs. This process will be all the more efficient than the graph has small blobs and long chains (hence it will be maximally efficient for trees). The second tool is a universal lower bound
\begin{equation}
\frac{\ell^D}{\Delta(G)}\le d^D\label{lowerbound1}
\end{equation}
where $\Delta(G)$ is the degree of $G$, i.e. the maximal degree of its vertices. It will be obtained by splitting the graph into the disconnected union of its edges and defining a new spectral triple on this split version. An embedding of the canonical spectral triple into the ``split'' one will yield the lower bound.

The paper is organized as follows. In section \ref{secrev} we briefly review the known properties of $d^D$ on graphs. In section \ref{secref} we introduce some refinements: first we give a rough estimate of the modification of $d^D$ introduced by the removal of an edge, and then we introduce the blob-chain decomposition and use it to estimate the noncommutative distance between two points in the graph. This raises the question of estimating the noncommutative length of chains. This we do in section \ref{secestim}, where we also prove \eqref{lowerbound1}. Finally in section \ref{conclusion} we will list some of the problems we see as the most important in view of further explorations.


All the graphs considered in this paper are connected, unless specified otherwise.

\section{Noncommutative distances in graphs: a short review}\label{secrev}

\subsection{General definitions}
We consider an undirected graph $G=(V,E)$ with no loops, where $V=\{1,\ldots,n\}$ is the set of vertices and $E$ is the set of edges. An edge is a pair $\{i,j\}$ with $i\not= j$. For some purposes it will be useful to introduce an orientation on $G$, which is a pair of functions $s,t$ (source, target), from $E$ to $V$, such that $\{s(i,j),t(i,j)\}=\{i,j\}$ for all edge $\{i,j\}$.  A \emph{subgraph} of $G$ is a graph $H=(V',E')$ such that $V'\subset V$ and $E'\subset E$. If $V'$ is a subset of $V$, the subgraph \emph{induced} by  $V'$  is the graph $(V',E')$ such that $E'$ contains each edge of $E$ with endpoints both  in $V'$. The \emph{degree} of a vertex is the number of edges incident on it. We will use the convention that a \emph{path} in a graph is sequence of edges which join a sequence of vertices \emph{with no   repetition of edge} (this is sometimes called a \emph{trail}). We will say that a path is \emph{simple} if there is no repetition of vertices. A path which is also an induced graph is called an \emph{induced path}. A graph which consists of a single simple path will be called a \emph{chain}.

We want to encode a graph in a spectral data. For this, we build from $V$  the $C^*$-algebra $\A=\CC^V$ and the Hilbert space $\H=\CC^V$. The algebra is represented on $\H$ by pointwise multiplication, i.e. we have a representation $\pi : \A\rightarrow \End(\H)\simeq M_n(\CC)$ such that $\pi(\A)$ is the algebra of diagonal matrices. So far we haven't made use of $E$. The set of edges plays the role of a differential structure on $V$ \cite{dimakis}. We thus define the bimodule of \emph{1-forms} $\Omega^1$ to be the set of matrices $\omega\in M_n(\CC)$ such that 
\begin{equation}
\{i,j\}\notin E\Rightarrow \omega_{i,j}=\omega_{j,i}=0
\end{equation}
(Note in particular that the diagonal elements of $\omega$ vanish.) The data $\B(G):=(\A,\H,\Omega^1)$ is the so-called \emph{algebraic background} which is canonically attached to the graph $G$. We will not make further use of this notion and we refer to \cite{proceedingcorfu} for more details. Nevertheless, we will want our Dirac operators to respect the differential structure incarnated by $\Omega^1$. First, let us recall that a Dirac operator is in our context simply a self-adjoint matrix $D\in \End(\H)$. We will say that such a $D$ is \emph{compatible} with $\Omega^1$ iff
\begin{equation}
\forall a\in \A,\ [D,\pi(a)]\in\Omega^1
\end{equation}
This is equivalent to require that the elements $D_{ij}$ satisfy:
\begin{equation}
(i\not= j\mbox{ and }\{i,j\}\notin E)\Rightarrow D_{i,j}=D_{j,i}=0
\end{equation}
Thus it is possible that $D$ has non-vanishing diagonal elements. A compatible Dirac operator  uniquely defines  a set of   weights on the edges of $G$: the weight  of the edge $\{i,j\}$ is by definition $w^D_{ij}:=|D_{ij}|^{-1}$ and belongs to $]0;+\infty]$. Out of the weights one can define the geodesic distance on $V$ which we will denote by $\ell^D$. By definition the length $\ell^D({\cal P})$ of a path in $G$ in the sum of the weights of its edges, and the geodesic distance $\ell^D(i,j)$ is the infimum of the lengths of all paths from $i$ to $j$. Note that $\ell^D$ satisfies all the axioms of a distance except for the fact that $\ell^D(i,j)$ may be infinite.

It turns out that we can define another distance on $V$, by the formula:

\begin{equation}
d^D(i,j)=\sup_{a\in\A}\{|a(i)-a(j)|, \|[D,\pi(a)]\|\le 1\}\label{formuledist}
\end{equation}

We call $d^D$   the \emph{noncommutative (NC) distance} defined by $D$. It is  a distance in the same generalized sense as $\ell^D$.  When $D$ is clear by the context, we write the NC and geodesic distances simply $d$ and $\ell$. An immediate computation shows that for a graph with two vertices, $d=\ell$. A general theme in what follows will be to compare $d$ and $\ell$ for more complicated graphs.

Let us end this section by observing that the space of Dirac operator compatible with $\Omega^1$ is larger than the space of weights on $G$, because of the diagonal entries and the phases. The former play no role in the computations of distances, so we will always suppose them to vanish. The latter do play a role, but we will often  have to  assume that the entries of $D$ are non-negative real numbers.

\subsection{General properties of $d$}\label{secreview}
This section borrows results from \cite{ikm} to which we refer for the proofs. We start with the following:
\begin{lemma}
In formula \eqref{formuledist} one can take $a$ to be real, positive, and such that $\|[D,\pi(a)]\|=1$.
\end{lemma}
We now come to distance computations \emph{per se}.
\begin{lemma}\label{cancellemma}
If $D'$ is obtained by setting to zero the $i$-th line and $i$-th column of $D$, then $d^{D'}\ge d^D$.
\end{lemma}
Observe that if $V'\subset V$, then the submatrix $D_{V'}$ of $D$ obtained by erasing all the lines and columns corresponding to vertices not in $V'$ is a Dirac operator on the subgraph induced by $V'$. Let us call $d^D_{V'}$ the associated Connes' distance. Repeated applications of lemma \ref{cancellemma} show that 
\begin{equation}
\forall i,j\in V',\ d^D_{V'}(i,j)\ge d^D(i,j)
\end{equation}
In particular we can take for $V'$ a two-elements set. We thus obtain    $d^D(i,j)\le w^D(i,j)$. Now since $d^D$ satisfies the triangle inequality, we have $d^D(i,j)\le \sum_{(i_k,j_k)\in{\cal P}} d^{D}(i_k,j_k)$, where the sum extends over a path ${\cal P}$ leading from $i$ to $j$ and $(i_k,j_k)$ are the edges of this path. Hence $d^D(i,j)\le\ell^D({\cal P})$ for any path joining $i$ and $j$. We therefore obtain the following important corollary:

\begin{theorem}
For any $D$ and any $i,j\in V$ one has
\begin{equation}
d^{D}(i,j)\le \ell^D(i,j)
\end{equation}
\end{theorem}

Let us now turn to the relations between $d^D$ and paths in $G$. The most basic is the following:

\begin{theorem}\label{connected}
$\forall i,j\in V$, $d^D(i,j)<\infty$ iff $i$ and $j$ are connected.
\end{theorem}

In order to state the next key property, let us introduce some notations. For any two $i,j\in V$, call $\PP(i,j)$ the subgraph of $G$ which is the union of all paths joining $i$ and $j$. 

\begin{theorem}\label{pathdecomp}
For any two vertices $i,j\in V$, the distance $d^D(i,j)$ is equal to the distance $d^{D'}(i,j)$ where $D'$ is the Dirac operator induced by $D$ on $\PP(i,j)$.
\end{theorem}

Let us give a first application of this theorem. Recall that an edge $(i,j)$ is a \emph{bridge} if removing it disconnects the graph. Since $(i,j)$ is then the  only path from $i$ to $j$ we obtain:

\begin{cor} Let $(i,j)$ be a bridge. Then $d^D(i,j)=w^D(i,j)=\ell^D(i,j)$.
\end{cor}

It turns out that theorem \ref{pathdecomp} is still true if we use simple paths instead of paths (but the proof in \cite{ikm} must be slightly modified), as we will see in section \ref{blobchaindec}.

\subsection{Specific cases}\label{secomput}
In this section we gather the explicit computations of $d$ we could find in the literature. The following result is proven in \cite{ikm}.


\begin{propo}\label{firstexamples}
\begin{enumerate}
\item If $n=3$ and $D$ has real entries, then
$$d^D(1,2)=\sqrt{\frac{D_{13}^2+D_{23}^2}{D_{12}^2D_{13}^2+D_{12}^2D_{23}^2+D_{23}^2D_{13}^2}}.$$
\item If $G$ is complete and the weights $w^D(i,j)$ are all equal, then
$$d^D(i,j)=w\sqrt{\frac{2}{n}}$$ 
where $w$ is the common weight.
\item If $G$ is complete and the weights are all equal to $w$ except $w^D(1,2)$ which is equal to $\infty$, then
$$d^D(1,2)=w\sqrt{\frac{2}{n-2}}$$
\end{enumerate}
\end{propo}

The following negative result is also worth mentioning (we give here a simplified statement,  for more details see \cite{ikm}).

\begin{propo}\label{negative}
Noncommutative distances in a 4-point graph are generically not computable : they are roots of polyomials not solvable by radicals. However, they \emph{are} computable when the graph is a cycle.
\end{propo}

In order to state the next results we need to introduce yet more terminology. Let $C$ be a chain. Its graph-theoretic length $|C|$ is just the number of its edges. Suppose $C$ has $n$ vertices   and is equipped with weights $w_1,w_2,\ldots,w_{n-1}$ on its edges. Such a weighted chain  will be denoted $C=w_1-w_2-\ldots-w_{n-1}$. There uniquely corresponds to it a Dirac operator $D$ with positive real entries  $D_{i,i+1}=D_{i+1,i}=w_i^{-1}$ for $i=1,\ldots,n-1$. We write
\begin{equation}
\lambda(C):=d^D(1,n)
\end{equation}
which we call the \emph{noncommutative length} of $C$. There are thus 3 distinct lengths for $C$ which should not be confused: $|C|$, $\lambda(C)$ and $\ell(C)$ which are in order the graph-theoretic, noncommutative and geodesic lengths.  

The noncommutative lengths of chains with $1$, $2$ and $3$ edges can be explicitly computed (for $|C|=2$ see proposition \ref{firstexamples} with $D_{12}=0$, for $|C|=3$, see appendix \ref{chain3}). For longer chains with generic weights, the explicit computation is not possible. However it can be performed for a chain with equal weights   \cite{bimonte1994distances,atzmon1996distances}. We gather the results below.

\begin{propo}
We have:
\begin{itemize}
\item $\lambda(w_1)=w_1$,
\item $\lambda(w_1-w_2)=\sqrt{w_1^2+w_2^2}$,
\item  $\lambda(w_1-w_2-w_3)=\frac{\sqrt{w_1^2+w_2^2}\sqrt{w_3^2+w_2^2}}{w_2}$, if  $w_2>\sqrt{w_1w_3}$, else $\lambda(w_1-w_2-w_3)=w_1+w_3$.
\end{itemize}
Moreover, if all the weights are equal to $w$, one has for all $k\in\NN^*$:
\begin{eqnarray}
\lambda(C)&=&w\sqrt{k(k+1)}\mbox{ when }|C|=2k\cr
\lambda(C)&=&wk\mbox{ when }|C|=2k-1.\label{longueurschainescstantes}
\end{eqnarray}
\end{propo}
Formulae \eqref{longueurschainescstantes} also apply to the noncommutative distance between two vertices of a rectangular lattice lying on the same line \cite{bimonte1994distances}.

\section{Refinements}\label{secref}
In this section we gather results which are minor improvements on those reviewed in section \ref{secreview}. The general idea is to treat, to some extent, the noncommutative distance as a modified geodesic distance. Let us start with an observation. From lemma \ref{cancellemma} we can infer the upper bound
\begin{equation}
d(x,y)\le \min_{G'}d'(x,y)
\end{equation}
where the min extends over all subgraph $G'$ of $G$ and $d'$ is computed by restricting $D$ to $G'$. In particular we have
\begin{equation}
d(x,y)\le \min_{\gamma : x\rightarrow y}\lambda(\gamma)\label{majinducedpaths}
\end{equation}
where $\gamma$ runs over all induced paths from $x$ to $y$ in $G$. This inequality is obviously an equality when $G$ is a tree. In subsection \ref{blobchaindec} we will use a certain decomposition of the graph to exploit this trivial fact as much as possible.

However, one must not push the analogy too far: the geodesic and noncommutative distances behave very differently. First, the inequality in \eqref{majinducedpaths} is generally strict, even for graphs as simple as a triangle. In the latter case, when all the weights are equal to $1$, we see from proposition \ref{firstexamples} that the LHS of \eqref{majinducedpaths} is $\sqrt{2/3}$ while the RHS is $1$.  We thus see that even though the direct link is, in the noncommutative sense, the shortest path between $x$ and $y$, the weights on the other edges still matter. Another example of the marked difference between $d$ and $\ell$ is provided by the noncommutative lengths of chains,  which  is not additive with respect to concatenation, as can be seen  from  \eqref{longueurschainescstantes}.

Still, the noncommutative and geodesic distances do share a property: they increase when one removes a vertex and all the edges attached to it. The geodesic distance $\ell(x,y)$ also increase when a single edge, say $(1,2)$  is removed. Would it be true also for the noncommutative distance ? The next subsection is devoted to this question.

\subsection{Adding and removing edges}

We can see from the explicit computation in  proposition \ref{firstexamples} that $d$ increases when we remove one edge in the triangle graph. The reason behind this fact is that the operator norm of a real $3\times 3$ antisymetric matrix is proportional to its $L^2$-norm, and the latter increases with the absolute values of the entries. Hence, if  $a\in \A$ satisfies $\|[D,a]\|=1$ and $|a(x)-a(y)|=d(x,y)$ then $\|[D',a]\|\le 1$ where $D'$ is $D$ with the $(1,2)$ and $(2,1)$ entries set to zero. Thus $d'(x,y)\ge d(x,y)$ where $d'$ is the NC distance on the weighted graph amputated from    $(1,2)$ edge. However this property of antisymmetric matrices does not extend to $n\ge 4$. Thus, it cannot be asserted in general that $\|[D,a]\|\le 1\Rightarrow\|[D',a]\|\le 1$.   Still, numerical simulations with hundreds of up to $5\times 5$ random Dirac operators seem to indicate that $d'(x,y)\ge d(x,y)$. The case of a complete graph in proposition \ref{firstexamples} and theorem \ref{connected} also go in this direction. Yet, we  have failed to prove this property   and so we  leave it as an   open problem. Instead we will content ourselves with a crude estimate showing, at least, that $d'(x,y)$ cannot be much smaller that $d(x,y)$. 

\begin{propo}\label{addcut}
Let $d$ and $d'$ be the noncommutative distances between two vertices $x,y$ in weighted graphs $G$ and $G'$ which differ only by the addition or deletion of the unconnected edges  $e_1=(i_1,j_1),\ldots,(i_k,j_k)$ of weights $w_1,\ldots,w_k$. Then we have
\begin{equation}
\frac{ d}{1+m}\le d'\le (1+m')d\label{estimcut}
\end{equation}
where $m=\max_{s=1,\ldots,k}\frac{d(i_s,j_s)}{w_s}$ and $m'=\max_{s=1,\ldots,k}\frac{d'(i_s,j_s)}{w_s}$. In the particular case where one adds/remove a single edge of weight $w$ which connects $x$ and $y$, one obtains
\begin{equation}
|\frac{1}{d}-\frac{1}{d'}|\le \frac{1}{w}
\end{equation}
\end{propo}
\begin{demo}

We let $E$ be the hermitian matrix such that $E_{12}=e^{i\theta_{12}}/w_{12}$, where $w_{12}>0$ and $\theta_{12}\in\RR$, $E_{21}=E_{12}^*$ and $E_{ij}=0$ for every other pair of indices.

\begin{itemize}
\item First case : addition of the edge $(1,2)$. Hence $D_{12}=0$ and $D'=D+E$. For every $a\in \A$, one has $\|[D',a]\|\le \|[D,a]\|+\|[E,a]\|=\|[D,a]\|+\frac{|a(1)-a(2)|}{w_{12}}$. Suppose $a$ is such that $\|[D,a]\|=1$ and $|a(1)-a(n)|=d(1,n)$. Then one also has $|a(1)-a(2)|\le d(1,2)$ and we obtain
\begin{equation}
\|[D',a']\|\le 1,\mbox{ with }a'=\frac{a}{1+\frac{d(1,2)}{w_{12}}}
\end{equation}
Hence $|a'(1)-a'(n)|\le d'(1,n)$ and we find
\begin{equation}
\frac{w_{12}d(1,n)}{w_{12}+d(1,2)}\le d'(1,n)\label{ddprim}
\end{equation}
Since we can write $D=D'-E$, we find in a symmetric way
\begin{equation}
\frac{w_{12}d'(1,n)}{w_{12}+d'(1,2)}\le d(1,n)\label{dprimd}
\end{equation}
\item Second case : deletion of the edge $(1,2)$. Here $D_{12}=e^{i\theta_{12}}/w_{12}$ and $D'=D-E$. Equations \eqref{ddprim} and  \eqref{dprimd} are proven in the same way.
\item Third case : we add/remove several  edges $e_1=(i_1,j_1),\ldots,(i_k,j_k)$ of weights $w_1,\ldots,w_k$ which are not connected to each other. Since the edges are unconnected, the perturbation matrix $E$ is a direct sum whose norm is $\max_{s=1,\ldots,k}\frac{|a(i_s)-a(j_s)|}{w_s}$. The same proof as above yields the result.
\end{itemize}
\end{demo}

\begin{figure}[hbtp]
$$\entrymodifiers={[o][F-]}\xymatrix{ & \bullet \ar@{-}[dr]^{1} & & & \bullet\ar@{-}[dr]^{30} & \cr
\star \ar@{-}[ur]^{30} \ar@{-}[dr]^{1}& & \bullet \ar@{-}[r]^{1} & \bullet \ar@{-}[ur]^{1}   \ar@{-}[dr]^{1} & & \star\cr
& \bullet\ar@{-}[ur]^{1} & & & \bullet\ar@{-}[ur]^{1}  & }$$
\caption{A weighted graph with two ``heavy'' edges unconnected to one another}\label{heavy}
\end{figure}
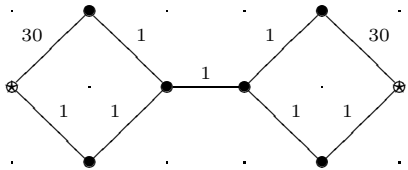

As an example of application of this proposition, let us consider the weighted graph of figure \ref{heavy}. The two stars will play the role of $x$ and $y$ in the proposition. The geodesic distance between the two ends of each edge of  weight $30$ is $3$ whether one removes these edges or not. Hence  one has $\max(m,m')\le 0.1$, thus one obtains $\frac{10}{11}d'\le d\le \frac{11}{10}d'$ from \eqref{estimcut}. Moreover, $d'$ can be computed explicitly from theorem \ref{pathdecomp} and formula \ref{longueurschainescstantes} (one finds $d'=3$).

%

\subsection{The blob-chain decomposition}\label{blobchaindec}
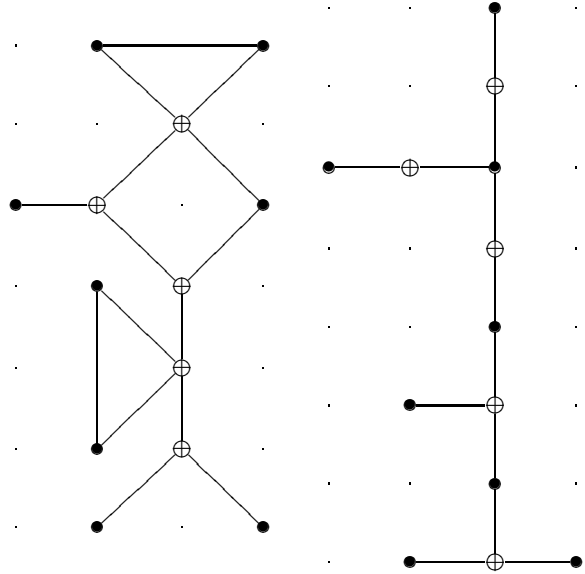
\begin{figure}[hbtp]
\begin{center}
\parbox{4cm}{$\entrymodifiers={[o][F-]}\xymatrix{ 
& \bullet\ar@{-}[rr]\ar@{-}[dr] &  &\bullet\ar@{-}[dl] \cr
 & & +\ar@{-}[dr] \ar@{-}[dl]   &\cr
 \bullet\ar@{-}[r] & +\ar@{-}[dr]  & &\bullet\ar@{-}[dl]\cr
 &\bullet\ar@{-}[dr] \ar@{-}[dd] & +\ar@{-}[d] & \cr
 & & + \ar@{-}[dl]\ar@{-}[d]&  \cr
  &\bullet &+ \ar@{-}[dl]\ar@{-}[dr] &  \cr
  & \bullet & & \bullet
}$}
\parbox{4cm}{$\entrymodifiers={[o][F-]}\xymatrix{ 
& & \bullet\ar@{-}[d]&\cr
& & +\ar@{-}[d]&\cr
\bullet\ar@{-}[r] & +\ar@{-}[r] & \bullet\ar@{-}[d]&\cr
& & +\ar@{-}[d]&\cr
& & \bullet\ar@{-}[d]&\cr
& \bullet\ar@{-}[r] & +\ar@{-}[d]&\cr
& & \bullet\ar@{-}[d]&\cr
&\bullet\ar@{-}[r] & + \ar@{-}[r]&\bullet\cr
}$}
\end{center}     
\caption{A graph and its block-cutpoint tree. The symbol $\oplus$ represents a cutpoint.} \label{bcdecomp}
\end{figure}

We are now going to introduce a decomposition of a graph  useful to estimate noncommutative distances.  For the graph theoretic notions in this section we refer to \cite{Harary}. A graph is called \emph{2-connected} if it is non-trivial and cannot be disconnected by the removal of a single vertex\footnote{Physicists call such graphs \emph{1-particle irreducible.}}. A \emph{block} of $G$ is a maximal 2-connected  subgraph of $G$. A \emph{cutpoint} of $G$ is a vertex $v$ such that the subgraph induced by $V\setminus\{v\}$ is disconnected. Every graph can be decomposed in blocks and cutpoints in such a way that two blocks cannot share more than one cutpoint. Note however, that a cutpoint can belong to several blocks (like the center of a star-shaped graph), and a block can contain several cutpoints. A block which has just one edge is a bridge. The structure of the graph $G$ can thus be summarized by a tree, called the \emph{block-cutpoint tree} $bc(G)$ defined in the following way:
\begin{itemize}
\item the vertices of $bcG)$ are the blocks and cutpoints of $G$,
\item each edge of $bc(G)$ joins a block $b$ with a cutpoint $c$, and does so iff $c\in b$.
\end{itemize}
Figure \ref{bcdecomp} exemplifies this construction. Note that there is a surjective map $bc$ from $V$ to the vertices of $bc(G)$ which sends a vertex $v$ to the block it belongs to if this block is unique, and to itself if $v$ is a cutpoint.

Our first result is a refinement of theorem \ref{pathdecomp} and its proof is entirely similar. 

\begin{definition}\label{pruning}
Let $i\not= j$ be two vertices in $G$. Let $\gamma$ be the set of vertices in the unique path from $bc(i)$ to $bc(j)$. Then the $(i,j)$-\emph{pruning} of $G$ is the subgraph $\Pi(i,j)$ of $G$ whose vertex set is the union of all the blocks and cutpoints in $\gamma$.
\end{definition}

Note that the $(i,j)$-pruning is the union of all simple paths from $i$ to $j$. There exist efficient algorithms to find the block-cutpoint tree (see for instance \cite{algo447}) which can also be used to prune a graph.


\begin{theorem}\label{pruneth}
For any two vertices $i,j\in V$, the distance $d^D(i,j)$ is equal to the distance $d^{D'}(i,j)$ where $D'$ is the Dirac operator induced by $D$ on the pruning $\Pi(i,j)$.
\end{theorem}
\begin{demo}
Let $V=V'\coprod V''$ where $V'$ is the set of vertices belonging to the pruning. For each $w\in V''$ there is a unique $w'\in V'$ with property that any path between $w$ and a $v\in V'$ goes through $w'$. For every function $b$ on $V'$, there is a unique extension $\tilde b$ on $V$ such that for all $w\in V''$, $\tilde{b}(w)=b(w')$. Now observe that in the block-decomposition of $\H$ induced by the partition $V=V'\coprod V''$, one has $[D,\tilde{b}]=[D',b]\oplus 0$. Hence, if we choose $b$ such that $\|[D',b]\|\le 1$ and $d^{D'}(i,j)=|b(i)-b(j)|$, then $\|[D,\tilde{b}\|\le 1$ hence $d^D(i,j)\ge |\tilde{b}(i)-\tilde{b}(j)|=|b(i)-b(j)|=d^{D'}(i,j)$. We conclude by lemma \ref{cancellemma}.
\end{demo}

In view of this result, we can now focus on evaluating distances $d^D(i,j)$ on graphs which have already been $(i,j)$-pruned. Hence we will make the following hypotheses:

\begin{itemize}
\item{H0} $i=1$ and $j=n$.
\item{H1} $1$ and $n$ belong to different blocks $b_1$ and $b_k$.
\item{H2} The graph $bc(G)$ is a chain with $b_1$ and $b_k$ as endpoints.
\end{itemize}

\begin{figure}[hbtp]
\begin{center}
$\entrymodifiers={[o][F-]}\xymatrix{ 
  & \bullet\ar@{-}[dr] &  & \bullet\ar@{-}[dr] &   &   &   &   & \bullet\ar@{-}[dr]\ar@{-}[dl] &   &   & \bullet\ar@{-}[dr]^>{\gamma_4}\ar@{-}[dl] &   \cr
\bullet\ar@{-}[ur]^<{\gamma_1}\ar@{-}[dr]\ar@{-}[r] & \bullet \ar@{-}[r]& \bullet\ar@{-}[dr]\ar@{-}[ur] & & \bullet\ar@{-}[r]^<{\gamma_1'} & \bullet \ar@{-}[r]\ar@{--}[u]& \bullet\ar@{-}[r]^>{\gamma_2}  & \bullet\ar@{-}[r]& \bullet\ar@{-}[r]\ar@{-}[u]\ar@{-}[d] & \bullet\ar@{-}[r]^<{\gamma_2'}^>{\gamma_3} & \bullet &  & \bullet \cr
  & \bullet\ar@{-}[ur] &  & \bullet\ar@{-}[ur] &   &   &   &   & \bullet\ar@{-}[ur]\ar@{-}[ul] &   &   & \bullet\ar@{-}[ur]\ar@{-}[ul] &   \cr
 \ar@{<->}[rrrr]|{b_1} &  &  &   &  \ar@{<->}[rrr]|{C_1}  &   &   &   \ar@{<->}[rr]|{b_2} &   &  \ar@{<->}[r]|{C_2}  &  \ar@{<->}[rr]|{b_3}  &   & 
}$
\end{center}
\caption{An example of the blob-chain decomposition. Here $\gamma_1=1$ and $\gamma_4=n$. The graph is $(1,n)$-pruned. The dotted line stands for a branch which has been cut during the pruning operation: it leads to parts of the graph which do not affect the value of $d(1,n)$.}\label{fignot}
\end{figure}
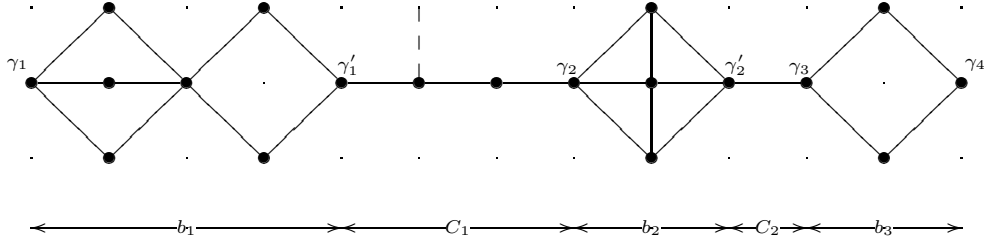

Note that some blocks may be connected by a  trivial chain, i.e. by a single cutpoint. Let us call a \emph{blob}   a maximal chain of such blocks, i.e. a maximal chain of blocks which does not contain any bridge. Let us now introduce a simplified version of the block-cut tree which is a single chain $C(G)$    defined by shrinking to a   vertex each blob. Note that the algebra of functions on $C(G)$ is isomorphic to the subalgebra of $\A$ containing the functions which are constant on blocks (hence on blobs, since the constant must be the same on two blocks which share a cutpoint). The chain $C(G)$ thus has the following structure:
$$C(G)= b_1-C_1-b_2-C_2-\ldots-C_{k-1}-b_k$$
where $b_1,\ldots,b_k$ are blobs and $C_1,\ldots,C_{k-1}$ are chains (this $k$ is not necessarily the same as before since there is in general less blobs than blocks, but we do not want to introduce a new notation). Observe that  if $k>1$ the blob $b_i$ is connected (in $G$) to $C_{i-1}$ by  a cutpoint $\gamma_{i}$ and if $k<n$ to $C_{i}$ by a cutpoint $\gamma_i'$. We complete this definition with $\gamma_1:=1$ and $\gamma_k':=n$. All these notations are summarized in figure \ref{fignot}.

Let us call for all $1\le i\le k$
\begin{equation}
\lambda(b_i)=\lambda(\gamma_{i},\gamma_{i}')
\end{equation}
the \emph{transverse length} of the blob $b_i$. Then a simple application of the triangle inequality  yields:
\begin{equation}
d(1,n)\le \sum_{i=1}^{k}\lambda(b_i)+\sum_{i=1}^{k-1}\lambda (C_i). \label{majorblobchain}
\end{equation}
We will obtain a minoration of $d(1,n)$ by considering a particular function $a$ on $G$ such that $\|[D,B]\|\le 1$. In order to define this function, we need one more piece of notation: for each chain $C_i$ in $C(G)$  we call $C_i'$ the chain $C_i$ amputated from its first and last edges. Hence if $C_1=w_1-w_2-\ldots-\ldots-w_{p-1}-w_p$, then $C_1'=w_2-\ldots-w_{p-1}$. If $\ell(C_i)\le 2$, the amputated chain $C_i'$ will be empty by definition. Let us now define $a$ inductively on each part of $C(G)$.
\begin{itemize}
\item \emph{Definition of $a$ on the blob $b_1$.} Let $B_1$ be the Dirac operator induced by $D$ on  $b_1$, and let $a_1$ be a function defined on the vertices of $b_1$ such that $\|[B_1,a_1]\|\le 1$ and $a_1(\gamma_1')-a_1(1)=d(1,\gamma_1')=\lambda(b_1)$. We let $a:=a_1$ on $b_1$.
\item \emph{Definition of $a$ on the chain $C_1$.} Recall first that $C_1$ starts at $\gamma_1'$ and ends at $\gamma_2$. Let $\Gamma_1'$ be the Dirac operator induced by $D$ on $C_1'$. Then we choose a function $a_1'$ on $C_1'$ such that $\|[\Gamma_1,a_1']\|\le 1$ and $a_1'($endpoint of $C_1')-a_1'($starting point of $C_1')=\lambda(C_1')$. We can suppose by adding a constant to $a_1'$ that its value on the starting point of $C_1'$ is equal to $a(\gamma_1')$. We let $a=a_1'$ on $C_1'$, and $a(\gamma_2)=a_1'($endpoint of $C_1')$. (Hence the value of $a$ remains constant on the two first (resp. two last) vertices of $C_1$.) We observe that $a(\gamma_2)-a(\gamma_1')=\lambda(C_1')$.
\item \emph{Definition of $a$ on $b_i$, $i>1$}. We note that $a$ is already defined on the entry point of $b_i$. Using the freeness of adding a constant we can extend it to   $b_i$ in such a way that $\|[B_i,a_{|b_i}]\|\le 1$,  and $a(\gamma_i')-a(\gamma_i)=\lambda(b_i)$.
\item \emph{Definition of $a$ on $C_i$, $i>1$}. Similarly, we can define $a_{|C_i}$ such that   $\|[\Gamma_i',a_{|C_i}]\|\le 1$,  and $a(\gamma_{i+1})-a(\gamma_i')=\lambda(C_i')$.
\end{itemize}
Now we have
\begin{eqnarray}
a(n)-a(1)&=&a(\gamma_k')-a(\gamma_k)+a(\gamma_k)-\ldots+a(\gamma_1')-a(\gamma_1)\cr
&=&\sum_{i=1}^k\lambda(b_i)+\sum_{i=1}^{k-1}\lambda(C_i').
\end{eqnarray}
Moreover, the Dirac operator $D$ can be written as
\begin{equation}
D=\begin{pmatrix}
B_1&{\tiny \begin{pmatrix}
0\cr
\vdots\cr
0\cr
w_1^{-1}
\end{pmatrix}}& &\cr
{\tiny\begin{pmatrix}
0&\ldots&0&w_1^{-1}
\end{pmatrix}} & \Gamma_1'&{\tiny\begin{pmatrix}
0\cr
\vdots\cr
0\cr
w_p^{-1}
\end{pmatrix}} & \cr
&{\tiny \begin{pmatrix}
0&\ldots&0&w_p^{-1}
\end{pmatrix}}& B_2&\cr
& & & \ddots
\end{pmatrix}
\end{equation}
where $w_1,\ldots,w_p$ are the weights of $C_1$. Since $a$ is constant on the first two and two last vertices of each chain, the commutator $[D,a]$ is the direct sum $[B_1,a_{|b_1}]\oplus 0\oplus [\Gamma_1',a_{|C_1'}]\oplus\ldots$, hence we see that $\|[D,a]\|\le 1$. We thus obtain
\begin{equation}
\sum_{i=1}^k\lambda(b_i)+\sum_{i=1}^{k-1}\lambda(C_i')\le d(1,n).\label{minorblobchain}
\end{equation}
Comparing \eqref{majorblobchain} and \eqref{minorblobchain} we see that the distance between to vertices in a graph is almost equal to the sum of the noncommutative lengths of the chains and blobs between them. The difference between this approximation and the exact result can be entirely attributed to the weights of the first and last edges of each chain. When these are large compared with the lengths of the blobs, another minoration might be useful:
\begin{equation}
\sum_{i=1}^{k-1}\lambda(C_i)\le d(1,n).
\end{equation}
It can be easily obtained using a function $a$ which is constant on each blob and realizes the length of each chain.  Combining these results we finally obtain the estimate
\begin{equation}
\max(\sum_{i=1}^{k}\lambda(b_i)+\sum_{i=1}^{k-1}\lambda(C_i'),\sum_{i=1}^{k-1}\lambda(C_i))\le  d(1,n)\le \sum_{i=1}^{k}\lambda(b_i)+\sum_{i=1}^{k-1}\lambda (C_i).
\end{equation}
In view of this result it is important to have estimate for  the NC length of chains. We will give one in the next  section.

\section{New estimates}\label{secestim}
In this section we prove two new estimates for the NC distance. One is completely general, the other only applies to chains. 

\subsection{A lower bound on the noncommutative distance}
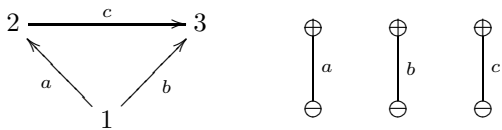
\begin{figure}[hbtp]
\begin{center}
\parbox{4cm}{$\xymatrix{ 
2\ar[rr]^c & & 3 \cr
 & 1\ar[ul]^a  \ar[ur]_b &
}$}\parbox{4cm}{$\entrymodifiers={[o][F-]}\xymatrix{+&+&+\cr
-\ar@{-}[u]_a&-\ar@{-}[u]_b&-\ar@{-}[u]_c }$}
\end{center}     
\caption{A directed graph   and its split version}\label{splitgraph}
\end{figure}

We already know an upper bound that applies to all weighted graphs, namely $d\le\ell$. In this section we  look for a lower bound of the same kind. For this we will use   a contruction that has been introduced in \cite{SST2} and which we recall here. The \emph{split graph} of $G$ is by definition the graph $\tilde G$ with vertex set $\tilde V:=E\times \{-,+\}$ and edges $\tilde{e}:=\{(e,-);(e,+)\}$ for all $e\in E$. It  can be seen as the disconnected sum of the edges of $G$ (see figure \ref{splitgraph}). However, in order to make this identification, we have to specify which endpoint of an edge $\tilde e$ in $\tilde G$ corresponds to which endpoint of $e$, and this amounts to fix an orientation of $G$. We thus consider $G$ to be oriented in this section. We denote by $e^-$ (resp. $e^+$) the source (resp. target) of the the edge $e$. There is thus a projection map $p : \tilde V\rightarrow V$ given by $(e,\pm)\mapsto e^\pm$. 
Let us now suppose that a compatible Dirac $D$ is given on $G$ and let us build a spectral triple over $\tilde G$. We first define the Hilbert space $\tilde\H=L^2(\tilde V)=\CC^E\otimes \CC^2$ which we equip with the canonical scalar product. The algebra $\A$ is the same as that used for $G$. However the representation now has multiplicities: the action of a function $a\in \A$ (which is defined on $V$) on an element $\Psi$ of $\tilde\H$ (which is a function on $\tilde V$) can be written as
\begin{equation}
\pi(a)\Psi(e,\pm):=a(e^\pm)\Psi(e,\pm).
\end{equation}
In other words $\pi(a)$ is the operator of multiplication by the pullback of $a$ by $p$. It will be useful to introduce yet another notation: for an edge $e$, we write $\H_{(e,\pm)}\subset \H$ for the subspace of functions vanishing outside $(e,\pm)$, $\H_e:=\H_{(e,-)}\oplus \H_{(e,+)}$ and for a vertex $v$, we $\H_v:=\bigoplus_{(e,\pm)|e^\pm=v}\H_{(e,\pm)}$. Using the identification $\H=\bigoplus_{e\in E}\H_e\simeq \bigoplus_{e\in E}\CC^2$, we can write $\pi(a)$ in matrix form as
\begin{equation}
\pi(a)=\bigoplus_e \begin{pmatrix}
a(e^-)&0\cr 0&a(e^+)
\end{pmatrix}\label{pia}
\end{equation}
whereas in the decomposition $\H=\bigoplus_{v\in V}\H_v$ we have
\begin{equation}
\pi(a)=\bigoplus_v a(v)I_{d(v)},\label{decoverv}
\end{equation}
where $d(v)$ is the degree of the vertex $v$, and also the multiplicity of the evaluation at $v$ as a subrepresentation of $\pi$.   Finally, the Dirac operator $\tilde D$ is
\begin{equation}
\tilde D=\bigoplus_{e\in E}\frac{1}{w^D(e)}\begin{pmatrix}0&1\cr 1&0\end{pmatrix}.
\end{equation}
The interest of this construction is that the pure state space of this new spectral triple is the same as before (since the algebra is the same), hence $\tilde D$ induces a distance on $V$, and it turns out that \cite{SST2}:
\begin{equation}
d^{\tilde D}=\ell^D.
\end{equation}
{\small For the sake of self-completeness, let us reproduce the proof here. One easily sees that 
  $\|[\tilde D,\pi(a)]\|=\sup_{e\in E}{|a(e^+)-a(e^-)|\over w_D(e)}$. If $\|[\tilde D,\pi(a)]\|\le 1$ and $(i_0,\ldots,i_k)$ is a geodesic path in $G$  joining $i_0$ and $i_k$, then we see that 

\bea
|a(i_k)-a(i_0)|&\le&\sum_{r=0}^{k-1}|a({i_{r+1}})-a({i_r})|\cr
&\le& \sum_{r=0}^{k-1}w_D(i_r,i_{r+1}) \cr
&=&\ell^D(i_0,i_k)
\eea
since $(i_0,\ldots,i_k)$ is a geodesic path. Moreover if we take $a$ to be the function $a(j):=\ell^D(i_0,j)$ then we have $a(i_k)-a(i_0)=\ell^D(i_0,i_k)$, so that we just need to prove that $\|[\tilde D,\pi(a)]\|\le 1$ for this particular function. This amounts to prove that for each edge $e=\{j_1,j_2\}$ one has $|\ell^D(i_0,j_1)-\ell^D(i_0,j_2)|\le w^D({j_1,j_2})$. But this is true because $|\ell^D(i_0,j_1)-\ell^D(i_0,j_2)|\le \ell^D(j_1,j_2)$ by the triangle inequality and $\ell^D(j_1,j_2)\le w^D({j_1,j_2})$.

}

We are now going to use this new triple to obtain a minoration of $d^D$. The idea is to observe that if we pull back functions on $V$ by the projection $p : \tilde V\rightarrow V$ we obtain a canonical embedding of $\H$ into $\tilde\H$. Let us call $\Sigma$ the image of this embedding: it contains so-called \emph{graph states} $\Psi$ such that $\Psi(x)=\Psi(y)$ whenever $p(x)=p(y)$. Then we consider the block-decomposition of $\tilde D$ with respect to $\Sigma$ and $\Sigma^\perp$.  For this, let us write $\delta_x$, $x\in\tilde V$ the canonical basis of $\tilde\H$. Then an orthonormal basis of $\Sigma$ is given by the following vectors indexed by $v\in V$:

\begin{equation}
e_v:=\frac{1}{\sqrt{\Delta(v)}}\sum_{x|p(x)=v}\delta_x
\end{equation}
where $\Delta(v)=|p^{-1}(v)|$ is the degree of $v$. Now we have
\begin{eqnarray}
\bra e_{v'},\tilde D e_{v}\ket&=&\frac{1}{\sqrt{\Delta(v)\Delta(v')}}\sum_{x,x'|p(x)=v,p(x')=v'}\bra \delta_{x'},D\delta_x\ket\cr
&=&\frac{1}{\sqrt{\Delta(v)\Delta(v')}} D_{v'v}
\end{eqnarray}
where we have used the fact that $\bra \delta_{x'},D\delta_x\ket$ vanishes unless $x$ and $x'$ are the two ends of a given edge, in which case this edge must link $v$ to $v'$. There are thus only two possibilities: $x=(e,+)$ and $x'=(e,-)$ or the opposite, only one of which is compatible with $p(x)=v$ and $p(x')=v'$.

Now let us complete $(e_v)$ to an orthonormal basis $B$ of $\tilde\H$. Calling $U$ the unitary  matrix gathering the elements of $B$, we have
\begin{equation}
\tilde D=U\begin{pmatrix}
\Delta^{-1/2} D \Delta^{-1/2} & ?\cr ?& ?
\end{pmatrix}U^\dagger
\end{equation} 
where $\Delta$ is the diagonal matrix containing the   degrees of the vertices. Moreover, if $a\in \A$, then clearly
\begin{equation}
\tilde\pi(a)e_v=a(v)e_v
\end{equation}
Since $\Sigma$ is stable by $\tilde\pi(a)$ and $\tilde\pi(a)^\dagger$, $\Sigma^\perp$ is also stable and we can conclude that
\begin{equation}
\tilde\pi(a)=U\begin{pmatrix}
\pi(a)&0\cr 0&?
\end{pmatrix}U^\dagger
\end{equation}
(A little more work shows that $?$ is also diagonal, but non-zero). Since unitary conjugation preserves the operator norm, we have
\begin{equation}
\|[\tilde D,\tilde\pi(a)]\|=\|\begin{pmatrix}
[D_\Delta,\pi(a)]&?\cr ?&?
\end{pmatrix}\|\ge \|[D_\Delta,\pi(a)]\|
\end{equation}
where $D_\Delta=\Delta^{-1/2}D\Delta^{-1/2}$. Hence $\{a|\|[\tilde D,\tilde\pi(a)]\|\le 1\}\subset \{a|\|[  D_\Delta,\pi(a)]\|\le 1\}$. Thus for all $i,j\in V$, $d^{\tilde D}(i,j)=\ell^D(i,j)\le d^{D_\Delta}(i,j)$.  In order to   go further, let $\Delta(G)$ be the maximal degree of the vertices of $G$, and observe that $\|\Delta^{-1/2}X\|\ge \Delta(G)^{-1/2}\|X\|$ for any vector $X$. Thus
\begin{eqnarray}
\|[D_\Delta,\pi(a)]\|&=&\|\Delta^{-1/2}[D,\pi(a)]\Delta^{-1/2}\|,\mbox{ since }[\Delta,\pi(a)]=0\cr
&=&\sup_{X\not=0} \frac{\|\Delta^{-1/2}[D,\pi(a)]\Delta^{-1/2}X\|}{\|\Delta^{-1/2}X\|}\frac{\|\Delta^{-1/2}X\|}{\|X\|}\cr
&\ge &\|\Delta^{-1/2}[D,\pi(a)]\|\Delta(G)^{-1/2}\cr
&\ge&\|[D,\pi(a)]\Delta^{-1/2}\|\Delta(G)^{-1/2},\mbox{ since }\Delta^{-1/2}\mbox{ and }i[D,\pi(a)]\mbox{ are selfadjoint}\cr
&\ge&\|[D,\pi(a)]\|\Delta(G)^{-1},\mbox{ as in the second step}
\end{eqnarray}
We therefore obtain the following result:
\begin{theorem}
Let $\Delta={\rm diag}(\Delta(1),\ldots,\Delta(n))$ be the diagonal matrix containing the degrees of the vertices,   $D^\Delta=\Delta^{-1/2}D\Delta^{-1/2}$ and $\Delta(G)=\max\{\Delta(i)|i=1,\ldots,n\}$. Then we have
\begin{equation}
\ell^D\le d^{D_\Delta}\le \Delta(G)d^D.\label{lowerbound}
\end{equation}
\end{theorem} 
 
\begin{rem} It is interesting to note that a lower bound similar to \eqref{lowerbound} has been found in \cite{requardt2002distance} for the noncommutative distance associated to a different spectral triple attached to the graph.
\end{rem}

\subsection{Estimating the noncommutative lengths of chains}
In view of section \ref{blobchaindec}, is particularly important to have good estimates for the noncommutative lengths of chains.

\begin{lemma}\label{lemma3} Let $b_1,\ldots,b_{n-1}\in\CC$ and $B\in M_n(\CC)$ be the  matrix such that $B_{i,i+1}=b_i$, $B_{i+1,i}=-\bar b_i$ for $i=1,\ldots,n-1$, and all other entries equal to $0$. Then the operator norm of $B$ satisfies the following inequalities

\begin{eqnarray}
\max_{1\le i\le n-2}\sqrt{|b_i|^2+|b_{i+1}|^2}&\le&\|B\|\cr
\frac{2}{n}\sum_{i=1}^{n-1}|b_i|&\le&\|B\|\cr
\|B\|&\le&\max_{1\le i\le n-2}\ |b_i|+|b_{i+1}| \cr
\|B\|&\le&2\cos(\frac{\pi}{n+1})\max_{1\le i\le n-1}|b_i|\label{opnormestim}
\end{eqnarray}
\end{lemma}
\begin{demo}
Let $X\in \CC^n$ such that $\|X\|=1$. Then 
\begin{equation}
\|BX\|^2=|b_1|^2|x_2|^2+|-\bar b_1x_1+b_2x_3|^2+\ldots+|-\bar b_{n-2}x_{n-2}+b_{n-1}x_n|^2+|b_{n-1}|^2|x_{n-1}|^2\label{BXsq}
\end{equation}
To obtain the first inequality we just choose $x_i=1$ for $i$ such that $|b_i|^2+|b_{i+1}|^2$ is maximal, and $x_j=0$ for $i\not=j$.

Let us now write $\tilde B$ for the real symmetric matrix obtained from $B$ by replacing each entry with its modulus. We claim that $\|B\|=\|\tilde B\|$. Indeed, we notice that the phases of $x_1,\ldots,x_{n}$ can be arranged such that the moduli in \eqref{BXsq} reach their maximum values, namely such that
\begin{equation}
\|BX\|^2=|b_1|^2|x_2|^2+(|b_1||x_1|+|b_2||x_2|)^2+\ldots\label{BXsq2}
\end{equation}
The RHS of \eqref{BXsq2} is clearly $\le \|\tilde B\|^2$. Hence we obtain $\|B\|\le\|\tilde B\|$. The other inequality is obtained by replacing $|x_1|,\ldots,|x_n|$ by the components of the normalized dominant eigenvector of $\tilde B$, which has non-negative entries by Perron-Frobenius theorem. For this vector $\|\tilde B\|$ is reached since $\tilde B$ is normal. To obtain the second inequality, let $U=\frac{1}{\sqrt{n}}(1,\ldots,1)^T$ and observe that $\bra U,\tilde BU\ket\le \|\tilde B\|$ by the min-max principle.

The third inequality is easily obtained  from   Gershgorin circles.

Finally the fourth inequality is   obtained by majorizing each $|b_i|$ and using the norm of the tridiagonal Toeplitz matrix $T=\begin{pmatrix}
0&1&& \cr
1&\ddots&\ddots& \cr
&\ddots&\ddots & 1\cr
 & & 1& 0
\end{pmatrix}$ which is $\cos\frac{\pi}{n+1}$ (\cite{gantmacherkrein}, p 111).
\end{demo}

 
Now we apply the lemma with $B=-[D,\pi(a)]$. This yields $b_i=\frac{X_i}{w_i}$, where $X_i=a_{i+1}-a_{i}$. The noncommutative length $\lambda(C)$  is the maximum of 
\begin{equation}
f(b_1,\ldots,b_{n-1})=\sum_{i=1}^{n-1}w_ib_i=a_n-a_1\label{funcf}
\end{equation}
subject to the constraint $\|B\|\le 1$.  Let us define:
\begin{eqnarray}
L_1(C)&=&\max \{f(b_1,\ldots,b_{n-1})||b_i|^2+|b_{i+1}|^2\le 1, \forall 1\le i\le n-2\}\cr
L_2(C)&=&\max \{f(b_1,\ldots,b_{n-1})|\sum_{i=1}^{n-1}|b_i|\le \frac{n}{2}\}\cr
R_1(C)&=&\max \{  f(b_1,\ldots,b_{n-1})||b_i|+|b_{i+1}|\le 1, \forall 1\le i\le n-2\}\cr
R_2(C)&=&\max \{  f(b_1,\ldots,b_{n-1})||b_i|\le\frac{1}{2\cos\frac{\pi}{n+1}}, \forall 1\le i\le n-1\}\cr
\end{eqnarray}
Note that here $b_i\in \RR$, but since $w_i\ge 0$ for all $i$, one can consider $b_i\ge 0$ in all the above expressions and suppress the absolute values. Clearly one has the estimate
\begin{equation}
\max(R_{1}(C),R_2(C))\le \lambda(C)\le \min(L_1(C),L_2(C))
\end{equation}
It is immediate to obtain
\begin{equation}
R_2(C)=\frac{\ell(C)}{2\cos\frac{\pi}{n+1}},\ L_2(C)=\frac{n}{2}\max_{1\le i\le n-1}w_i
\end{equation}
%
To obtain the values of $L_1$ and $R_1$ needs more work. This is done in appendix \ref{appL1}. We find:
\begin{eqnarray}
L_1&=&\max_{C=C_1-\ldots-C_{p+1}}\sum_{i=1}^{p+1}\sqrt{\Sigma_{\rm even}(C_j)^2+\Sigma_{\rm odd}(C_j)^2}\cr
R_1&=&\max_{C=C_1-\ldots-C_{p+1}}\sum_{i=1}^{p+1}\max(\Sigma_{\rm even}(C_j),\Sigma_{\rm odd}(C_j))
\end{eqnarray}
where $\Sigma_{\rm even/odd}(C_j)$ is the sum of weights of even/odd indices contained in the subchain $C_j$, and the maximum extends over admissible decompositions of $C$ into subchains. The admissibility conditions are detailed in the appendix. It is worth noticing that $R_1\ge R_\emptyset:= \max(\Sigma_{\rm even}(C),\Sigma_{\rm odd}(C))$. Which functions among $L_1,R_1,L_2,R_2$ yield the best estimate depends on the   assignments of weights. Table \ref{examples} displays some examples. 
%

\begin{table}[hbtp]
\begin{tabular}{|c|c|c|c|c|c|c|}
\hline 
Weights &   $R_1$ & $R_2$ & $\lambda$ & $L_1$ & $L_2$ \\ 
\hline 
\small{$1-1-1-\ldots$ $2k-1$ times }&  $k$ & $\frac{k-1/2}{\cos\frac{\pi}{2k+1}}$ & $k$ & \small{$\sqrt{k^2+(k-1)^2}$} & $k$ \\ 
\hline 
\small{$1-1-1-\ldots$ $2k$ times }&   $k$ & $\frac{k}{\cos\frac{\pi}{2k+2}}$ & \small{$\sqrt{k(k+1)}$} & $k\sqrt{2}$  &  $k+1/2$ \\ 
\hline 
\small{$2-1-2-1-2$} &  6 & $\approx 4.4$ & 6 & $\approx 6.3$ & 6 \\ 
\hline 
\small{$1-2-1-2-1$} &  4 & $\approx 3.9$ & $\approx 4.4$   & $\approx 5.1$ & 6 \\  
\hline 
\end{tabular}
\caption{On these examples we see that the best lower (resp. upper) bounds can be given by $R_1,R_2$ (resp. $L_1,L_2$) according to the case. The NC length $\lambda$ has been computed numerically in the last line.}\label{examples}
\end{table}
Note that the third example in table \ref{examples} can be easily generalized: if $w_1\ge w_2$ the NC length of the chain $(w_1-w_2)^k-w_1$ is $(k+1)w_1$ since in this case $R_1$ and $L_2$ coincide. It gives back formula \eqref{longueurschainescstantes} when $w_1=w_2$.

\section{Conclusion, outlook}\label{conclusion}
The methods and results layed out in this paper make it clear that controlling the non-locality of the noncommutative distance is a delicate matter: if a graph consists of single blob, we know essentially nothing beyond the general estimate $\frac{\ell}{\Delta(G)}\le d\le \ell$. Important examples of this kind are regular lattices. It seems clear that a vertex $z$ which is geodesically very far from $x$ and $y$ should not matter much, if at all, for the computation of $d(x,y)$. This is certainly a direction which is worth exploring. Another one would be to improve  proposition \ref{addcut} which does not look like the end of the story concerning the addition or deletion of edges. These are two of many lines of investigation we hope our work could serve as a motivation for.

\appendix
\section{Computation of $L_1$ and $R_1$}\label{appL1}
We recall that we are looking for the maximum of 
\begin{equation}
f(x_1,\ldots,x_{n-1})=\sum_{i=1}^{n-1}w_ix_i
\end{equation}
over the convex set $K$ defined by
\begin{eqnarray}
x_1^2+x_2^2&\le& 1\cr
x_2^2+x_3^2&\le& 1\cr
\vdots&&\cr
x_{n-2}^2+x_{n-1}^2&\le&1\label{defK}
\end{eqnarray}
The maximum must be reached on an extreme point of $K$. 
Let $I$ be a subset of $\{2,\ldots,n-2\}$ such that no two elements of $I$ are consecutive and call $K_I$ the subset of $K$ defined by 
\begin{eqnarray}
x_i^2+x_{i+1}^2&=&1\mbox{ for }i\notin I,\cr
x_i^2+x_{i+1}^2&<&1\mbox{ for }i\in I\label{defKI}
\end{eqnarray} 
and by $\bar K_I$ the closure of $K_I$ in the relative topology of $K$ (defined by $\le$ instead of $<$ in \eqref{defKI}). We note that $I\subset J\Rightarrow \bar K_I\subset \bar K_J$.

\begin{lemma}\label{extremepoints}
The set of extreme points of $K$ is $\bigcup_I K_I$ where $I$ runs over the  subsets of $\{2,\ldots,n-2\}$ with no consecutive  elements.
\end{lemma}
\begin{demo}
Let us first check that the conditions are necessary. Suppose $u=(u_1,\ldots,u_{n-1})\in K$ is such that $u_1^2+u_2^2<1$. Then let $u'=(k,u_2,\ldots,u_{n-1})$ and $u''=(-k,u_2,\ldots,u_{n-1})$ where $k=\sqrt{1-u_2^2}$. It is clear that $u',u''\in K$ and $u\in[u',u'']$. For similar reasons, the last inequality must also be an equality. Now suppose $u_i^2+u_{i+1}^2<1$ and let $P$ be the $2$-plane defined by $x_1=u_1,\ldots x_{i-1}=u_{i-1},\ldots,x_{i+1}=u_{i+1},\ldots,x_{n-1}=u_{n-1}$. Then $u$ must be an extreme point of $K\cap P$. Now $K\cap P$ is the intersection of the unit disk with the rectangle defined by the inequations $x_i^2\le 1-u_{i-1}^2$ and $x_{i+1}^2\le 1-u_{i+2}^2$. By hypothesis, $u$ lies in the interior of the disk, which implies that the rectangle is inside the disk and $u$ is one of its vertices. Thus $u$ satisfies $u_i^2=1-u_{i-1}^2$ and $u_{i+1}^2= 1-u_{i+2}^2$.

Let us now  show that the elements of $K_I$ are extreme.  Let us call $\pi_i$ the projection onto the $(x_i,x_{i+1})$-plane. If $i\notin I$, $\pi_i(K_I)$ is a subset of the unit circle, hence its elements are necessarily extreme points of $\pi_i(K)$ which is a subset of the unit disk. Suppose $x=tx'+(1-t)x''$ for $i\in]0,1[$, $x\in K_I$, $x',x''\in K$. Then we have $\pi_i(x')=\pi_i(x'')$ for all $i\notin I$. By the defining conditions of $I$, this means that $x_j'=x_j''$ for all $j=1,\ldots,n-1$. Hence $x'=x''$ and $x$ is extreme.
\end{demo}

Let us first compute the maximum of $f$ on $K_\emptyset$. This set can be parametrized by:
\begin{eqnarray}
x_1=x_3=\ldots&=&\sin\theta_1\cr
x_2=x_4=\ldots&=&\cos\theta_1
\end{eqnarray}
Expressed as a function of $\theta_1$, $f$ becomes $\tilde f(\theta_1)=\Sigma_{\rm even}(C)\cos\theta_1+\Sigma_{\rm odd}(C)\sin\theta_1$, where for any chain $C$ of weights $w_1,\ldots,w_\ell$ we define 
\begin{eqnarray}
\Sigma_{\rm even}(C)&=&w_2+w_4+\ldots\cr
\Sigma_{\rm odd}(C)&=&w_1+w_3+\ldots
\end{eqnarray}
The maximum of $\tilde f$ on $K_\emptyset$ is then easily computed to happen when  
\begin{eqnarray}
\sin\theta_1&=&\frac{\Sigma_{\rm odd}(C)}{\sqrt{\Sigma_{\rm odd}(C)^2+\Sigma_{\rm even}(C)^2}}\cr
\cos\theta_1&=&\frac{\Sigma_{\rm even}(C)}{\sqrt{\Sigma_{\rm odd}(C)^2+\Sigma_{\rm even}(C)^2}} 
\end{eqnarray}
Thus we have
\begin{equation}
\max_{K_\emptyset}f=\sqrt{\Sigma_{\rm odd}(C)^2+\Sigma_{\rm even}(C)^2}=\sqrt{(w_1+w_3+\ldots)^2+(w_2+w_4+\ldots)^2}
\end{equation}
\begin{definition} The \emph{bare length} of the chain $C$ is  
\begin{equation}
L_\emptyset(C):=\sqrt{\Sigma_{\rm odd}(C)^2+\Sigma_{\rm even}(C)^2}.
\end{equation}
We say that $C$ is \emph{extremal} if $L_1(C)=L_\emptyset(C)$.
\end{definition}
\begin{rem} It is intringuing that the bare length coincides with the length computed thanks to a different Dirac operator in \cite{requardt2002distance}.
\end{rem}

Note that an immediate application of the subadditivity of the Euclidean norm shows that the bare length is subadditive with respect to concatenations:
\begin{equation}
L_\emptyset(C-C')\le L_\emptyset(C)+L_\emptyset(C')
\end{equation}
Note that $L_1$ is also subadditive since the maximum of $f$ goes up when the ``connecting condition'' $x_{j}^2+x_{j+1}^2\le 1$ is removed, and that without this condition $\sum_{i=1}^jw_ix_i$ and $\sum_{i=j+1}^{n-1}w_ix_i$ can be maximized separately. 

Clearly, any chain of length $\ell\le 3$ is extremal. For $\ell=1$ (i.e. for a bridge), the geodesic, noncommutative and bare lengths coincide. For $\ell=2$ we recover the case of a triangle with one edge removed (set $D_{23}=0$ in the first part of proposition \ref{firstexamples}).

Now suppose $\ell>3$ and $C$ is not extremal. Then $f$ reaches its maximum on $K_I$ for some subset   $I\not=\emptyset$. If $I=\{i_1,\ldots,i_p\}$ let us define a decomposition of $C$ into subchains $C_1,\ldots,C_{p+1}$ so that   $C_1$ will have weights $w_1,\ldots,w_{i_1}$ and be of length $\ell_1=i_1$, $C_2$ will have weights $w_{i_1+1},\ldots,w_{i_2}$ and be of length $\ell_2=i_2-i_1$ and so on, and $C_{p+1}$ have weights $w_{i_p+1},\ldots,w_{n-1}$. Note that each subchain $C_i$ is of length at least $2$. The points in $K_I$ are parametrized by $p+1$ angles $\theta_1,\ldots,\theta_{p+1}$  defined such that 
\begin{eqnarray}
x_{2k+1}&=&\sin\theta_1\mbox{ for all }k\mbox{ s.t. }1\le 2k+1\le i_1\cr
x_{2k}&=&\cos\theta_1\mbox{ for all }k\mbox{ s.t. }1\le  2k\le i_1\cr
x_{2k+1}&=&\sin\theta_2\mbox{ for all }k\mbox{ s.t. }i_1+1\le 2k+1\le i_2\cr
x_{2k}&=&\cos\theta_2\mbox{ for all }k\mbox{ s.t. }i_1+1\le  2k\le i_2\cr
\vdots& &\vdots
\end{eqnarray}
submitted to the $p$ conditions
\begin{equation}
x_{i_j+1}^2+x_{i_j}^2<1,\ j=1,\dots,p
\end{equation}
which translated in terms of the $\theta$-variables yields
\begin{eqnarray}
\sin^2\theta_{j+1}&<&\sin^2\theta_{j}\mbox{ if }i_j\mbox{ is even,}\cr
\sin^2\theta_{j+1}&>&\sin^2\theta_{j},\mbox{ if }i_j\mbox{ is odd}\label{condsincos}
\end{eqnarray}
The function to maximize is now
\begin{equation}
\tilde f(\theta_1,\ldots,\theta_{p+1})=\sum_{j=1}^{p+1} (\tilde \Sigma_{\rm odd}(C_j)\sin\theta_j+\tilde \Sigma_{\rm even}(C_j)\cos\theta_j),\label{fdetheta}
\end{equation}
where the tilde over $S$ means that the parity of a weight is defined with respect to the total chain $C$. If the maximum of $f$ on the compact set $\bar K_I$ is reached in $K_I$ the gradient  of $\tilde f$ must vanish in the open set defined by \eqref{condsincos}. This means that each summand of \eqref{fdetheta} is maximal, and by the analysis of $K_\emptyset$ above, we obtain that
\begin{eqnarray}
\sin\theta_j&=&\frac{\tilde \Sigma_{\rm odd}(C_j)}{\sqrt{\Sigma_{\rm odd}(C_j)^2+\Sigma_{\rm even}(C_j)^2}}\cr
\cos\theta_j&=&\frac{\tilde \Sigma_{\rm even}(C_j)}{\sqrt{\Sigma_{\rm odd}(C_j)^2+\Sigma_{\rm even}(C_j)^2}} 
\end{eqnarray} 
Hence we obtain that if $f$ reaches its maximum in $K_I$ then it is necessary that the following  conditions hold for $j=1,\ldots,p$:
\begin{eqnarray*}
\frac{\tilde \Sigma_{\rm odd}(C_{j+1})^2}{\Sigma_{\rm odd}(C_{j+1})^2+\Sigma_{\rm even}(C_{j+1})^2}&>&\frac{\tilde \Sigma_{\rm odd}(C_j)^2}{\Sigma_{\rm odd}(C_j)^2+\Sigma_{\rm even}(C_j)^2},\mbox{ for odd }i_j,\cr
\frac{\tilde \Sigma_{\rm odd}(C_{j+1})^2}{\Sigma_{\rm odd}(C_{j+1})^2+\Sigma_{\rm even}(C_{j+1})^2}&<&\frac{\tilde \Sigma_{\rm odd}(C_j)^2}{\Sigma_{\rm odd}(C_j)^2+\Sigma_{\rm even}(C_j)^2},\mbox{ for even }i_j
\end{eqnarray*}
or equivalently
\begin{equation}
 \Delta_j<0\mbox{ for all even }i_j,\mbox{ and }  \Delta_j>0\mbox{ for all odd }i_j\label{discon1}
\end{equation}
where $\Delta_j$ is the determinant
\begin{equation}
\Delta_j:=\begin{vmatrix}
\tilde \Sigma_{\rm even}(C_j)&\tilde \Sigma_{\rm even}(C_{j+1})\cr \tilde \Sigma_{\rm odd}(C_j)&\tilde \Sigma_{\rm odd}(C_{j+1})
\end{vmatrix}
\end{equation}
We will say that the decomposition of $C$ induced by $I$ is \emph{admissible} when  \eqref{discon1} is satisfied. In this case, the maximum of $f$ on $K_I$ is  equal to the sum of the bare lengths of the subchains $C_1,\ldots,C_{p+1}$, hence we have
\begin{equation}
L_1(C)\ge \sum_{j=1}^{p+1}L_\emptyset(C_j)
\end{equation}
Note that $\sum_{j=1}^{p+1}L_\emptyset(C_j)>L_\emptyset(C)$ since an equality would imply that all the $\Delta_j$ vanish. It follows that $C$ is not extremal.

Now suppose the   conditions  \eqref{discon1} do not hold. Then  the maximum of $f$ on $K_I$ is reached on $\bar K_I\setminus K_I\subset\bar K_{I'}$ for some $I'\subsetneq I$, so that $\max_{K_I}f=\max_{K_{I'}}f$. Since $f$ must reach its maximum on one of the $K_I$ we obtain the following result.

\begin{propo}\label{thchains}
We have
\begin{equation}
L_1(C)=\max_{C=C_1\cup \ldots C_{p+1}}\sum_{j=1}^{p+1}L_\emptyset(C_j)
\end{equation}
where the maximum extends over all admissible decompositions of $C$ into subchains of length at least 2. In particular $C$ is extremal iff $C=C$ is the only such decomposition.
\end{propo}

Note that conditions \eqref{discon1} are invariant with respect to the reversal of the chain. Indeed, let us call $C'$ the chain $w_n-\ldots-w_1$. The subset $I$ defines a decomposition of $C'$ such that $i_j'=n-i_j$.  The condition at $i_j'$ for $C'$ involves the determinant
$$\Delta_j'=\begin{vmatrix}
\tilde \Sigma_{\rm even}'(C_{j+1})&\tilde \Sigma_{\rm even}'(C_j)\cr \tilde \Sigma_{\rm odd}'(C_{j+1})&\tilde \Sigma_{\rm odd}'(C_j)
\end{vmatrix}$$
Now if $n$ is odd, $\tilde S'_{\rm even/odd}=\tilde \Sigma_{\rm even/odd}$, so that $\Delta_j'=-\Delta_j$. But since the parities of $i_j'$ and $i_j$ are opposite to each other, the conditions are equivalent for $C$ and $C'$. On the other hand if $n$ is even, $\tilde S'_{\rm even/odd}=\tilde S_{\rm odd/even}$ so that $\Delta_j=\Delta_j'$, and the conditions are again equivalent.

Note also that proposition \ref{thchains} does not mean that there exists a decomposition of any given chain into extremal ones. For instance, consider the chain $C=2-1-1-2-1-4-2$. Then it admits a unique decomposition which is $(2-1-1-2)-(1-4-2)$. It can be seen that $2-1-1-2$ can be further decomposed into $(2-1)-(1-2)$. Despite that fact, $C$ cannot be decomposed in 3 parts. The problem is that the admissibility conditions are not compatible with refinements of decompositions. However, they are compatible with some form of coarse-graining, as the following lemma shows.

\begin{lemma}\label{coarsegraining} If there exists an admissible decomposition of $C$ into $k$ subchains, with $k>2$, then there exists an admissible decomposition of $C$ into $k-1$ subchains.
\end{lemma}
\begin{demo}
Let $C=C_1-\ldots-C_k$ be an admissible decomposition. For all $i,j$ we define \begin{equation}
\Delta_{ij}:=\begin{vmatrix}
\tilde \Sigma_{\rm even}(C_i)&\tilde \Sigma_{\rm even}(C_{j})\cr \tilde \Sigma_{\rm odd}(C_i)&\tilde \Sigma_{\rm odd}(C_{j})
\end{vmatrix}.
\end{equation}
For any two non-zero real numbers $x,x'$ we   write $x\sim x'$ iff $xx'>0$. Since the entries of $\Delta_{ij}$ are positive, it is easy to show that $\Delta_{ij}\sim \Delta_{jm}\Rightarrow \Delta_{ij}\sim\Delta_{jm}\sim \Delta_{im}$.

Let us assume that no decomposition of $C$ into $k-1$ subchains is admissible. We are going to prove that  $|C_2|,\ldots,|C_{k-1}|$ are all odd and $\Delta_{j-1,j}\sim\Delta_{j-1,j+1}$ for all $2\le j\le k-1$.

Suppose $|C_2|$ is even. Since $C_1-C_2-\ldots-C_k$ is admissible and $|C_1|$ has the same parity as $|C_1|+|C_2|$, it follows that $\Delta_{12}\sim\Delta_{23}$. Hence $\Delta_{12}\sim\Delta_{13}$, and $\Delta_{12}+\Delta_{13}\sim\Delta_{23}$. This shows that $(C_1-C_2)-C_3-\ldots$ is admissible, which contradicts the hypothesis. Thus $|C_2|$ is odd. It follows that  $|C_1|$ and $|C_1|+|C_2|$ have opposite parities, so that the admissibility conditions entail $\Delta_{12}\sim -\Delta_{23}$. Suppose $\Delta_{23}\sim \Delta_{13}$. Then $\Delta_{13}+\Delta_{23}\sim \Delta_{23}$ and this shows that $(C_1-C_2)-C_3-\ldots$ is admissible. Since this is excluded, we have $\Delta_{23}\sim -\Delta_{13}$ and consequently $\Delta_{12}\sim \Delta_{13}$.

Suppose now that $|C_2|,\ldots,|C_l|$ are all odd and $\Delta_{j-1,j}\sim\Delta_{j-1,j+1}$ for all $2\le j\le l$ and some $l\le k-2$. By assumption, the decomposition $C_1-\ldots-(C_l-C_{l+1})-\ldots$ is not admissible. This means that at least one of the two following conditions is \emph{not} satisfied:
\begin{eqnarray}
\Delta_{l-1,l}+\Delta_{l-1,l+1}&\sim&\Delta_{l-1,l}\cr
\Delta_{l,l+2}+\Delta_{l+1,l+2}&\sim&\Delta_{l+1,l+2}
\end{eqnarray}
However the first of these conditions is satisfied by hypothesis. We conclude that $\Delta_{l,l+2}\sim\Delta_{l+2,l+1}$. But this entails $\Delta_{l,l+1}\sim \Delta_{l+2,l+1}$. This means that the admissibility conditions of $C_1-\ldots-C_k$ between $C_l$ and $C_{l+1}$ on one hand and between $C_{l+1}$ and $C_{l+2}$ on the other are opposite, which is the case iff $|C_{l+1}|$ is odd. Moreover we have also shown that  $\Delta_{l,l+1}\sim\Delta_{l,l+2}$. We have hence shown by induction that $|C_2|,\ldots,|C_{k-1}|$ are all odd and $\Delta_{j-1,j}\sim\Delta_{j-1,j+1}$ for all $2\le j\le k-1$.

We thus have $\Delta_{k-2,k-1}\sim\Delta_{k-2,k}$. Hence $\Delta_{k-2,k-1}\sim \Delta_{k-2,k-1}+\Delta_{k-2,k}$ which is the only condition required for the admissibility of $C_1-\ldots-C_{k-2}-(C_{k-1}-C_k)$. This concludes the proof by contradiction.
\end{demo}

The following proposition is immediate from the lemma.

\begin{propo}\label{CNSextremal} $C$ is extremal iff there is no admissible decomposition of $C$ into two subchains.
\end{propo}

Observe that if the admissibility  conditions \eqref{discon1} hold,    even with large inequalities instead of strict ones, then one has $L_1(C)\ge \sum_{j=1}^{p+1}L_\emptyset (C_j)$. Hence if the chains $C_1,\ldots,C_{p+1}$ are extremal, we obtain $L_1(C)\ge \sum_{j=1}^{p+1}L_1(C_j)$. Since the converse   holds by subadditivity, we obtain the following useful result.

\begin{propo} Let $C=C_1-\ldots-C_{p+1}$ be a decomposition of $C$ into extremal subchains. If the \emph{gluing conditions}
\begin{equation}
\Delta_j\ge 0\mbox{ for all  }i_j\mbox{ odd, and }\Delta_j\le 0\mbox{ for all  }i_j\mbox{ even}\label{gluingcond}
\end{equation}
hold, then $L_1(C)=\sum_{i=1}^{p+1}L_1(C_i)=\sum_{i=1}^{p+1}L_\emptyset(C_i)$.
\end{propo}

Let us now consider the   case where all the weights are equal.

\begin{propo} A chain whose weights are all equal is extremal.
\end{propo}
\begin{demo}
We can suppose without loss of generality that the weights are all equal to 1. We write $C_p$ for the chain $1-1-\ldots-1$ ($p$ terms). If $p=2k$ then the decomposition $C_{2k}=1-1/1-1/\ldots/1-1$ clearly satisfies the gluing conditions. It follows that $\lambda(C_{2k})=k\sqrt{2}=\sqrt{k^2+k^2}$. Hence $C_{2k}$ is extremal. If $p=2k+1$, consider a 2-decomposition $C_{2k+1}=C_1-C_2$. If $|C_1|=2l$ with $l\in\NN^*$, then the admissibility condition is $\begin{vmatrix}
l& k-l\cr
l & k-l+1
\end{vmatrix}<0$ which is false. If $|C_1|=2l+1$ the condition reads  $\begin{vmatrix}
l& k-l\cr
l+1 & k-l
\end{vmatrix}>0$ which is also false. Hence no 2-decomposition of $C_{2k+1}$ is admissible, and we conclude by proposition \ref{CNSextremal}.
\end{demo}

Proposition \ref{thchains} and lemma \ref{coarsegraining} can be easily converted into an algorithm to compute the noncommutative length of $L_1$ for any given chain: order all decompositions of the chain by inclusion of the sets $I$ of decomposition points. Starting with the larger sets (the more refined decompositions), check the admissibility condition. If a decomposition $C_1-\ldots-C_k$ is admissible, store the number $\sum_{i=1}^kL_\emptyset(C_i)$. If at least one $k$-decomposition is admissible, check $k-1$-decompositions and go on, else stop  and return the max of the stored values.

Let us now look at $R_1$. This is the maximum of $f$ subject to
\begin{eqnarray}
x_1+x_2&\le& 1\cr
x_2+x_3&\le& 1\cr
\vdots&&\cr
x_{n-2}+x_{n-1}&\le&1\label{defKprim}
\end{eqnarray}
and we can add $x_i\ge 0$ for all $i$ since the weights $w_i$ are positive. It is thus a linear programming problem in standard form. The value of $R_2$ is thus the maximum of $f$ over all vertices of the compact convex polytope defined by \eqref{defKprim}. We follow the same method as for the computation of $L_1$. The only differences are the following:
\begin{itemize}
\item In lemma \ref{extremepoints}, we only obtain the inclusion of the set of extreme points into $\bigcup_I K_I$, but this is sufficient for the rest of the computation.
\item Conditions \eqref{discon1} are replaced with:
\begin{eqnarray}
\tilde\Sigma_{\rm odd}(C_j)<\tilde\Sigma_{\rm even}(C_j),&\tilde\Sigma_{\rm even}(C_{j+1})<\tilde\Sigma_{\rm odd}(C_{j+1}),\mbox{ for }i_j\mbox{ even},\cr
\tilde\Sigma_{\rm odd}(C_j)>\tilde\Sigma_{\rm even}(C_j),&\tilde\Sigma_{\rm even}(C_{j+1})>\tilde\Sigma_{\rm odd}(C_{j+1}),\mbox{ for }i_j\mbox{ odd}\label{condR1}
\end{eqnarray}
\item $L_\emptyset$ is replaced with $R_\emptyset=\max(\Sigma_{\rm odd},\Sigma_{\rm even})$.
\end{itemize}

\section{Chains of length 3}\label{chain3}
In this section we explicitly compute the noncommutative length of the chain $C=w_1-w_2-w_3$. From  \eqref{funcf}, it is the maximum of $f(b_1,b_2,b_3)=w_1b_1+w_2b_2+w_3b_3$ submitted to the condition $\|B\|\le 1$ with 
\begin{equation}
B=\begin{pmatrix}
0&b_1&0&0\cr
-b_1&0&b_2&0\cr
0&-b_2&0&b_3\cr
0&0&-b_3&0
\end{pmatrix}
\end{equation}
The characteristic equation of this matrix is of the form $P(\lambda^2)=0$, so that the norm of $B$ can be easily computed. We find
\begin{equation}
\|B\|=\sqrt{\frac{\beta+\sqrt{\beta^2-4b_1^2b_3^2}}{2}}
\end{equation}
with $\beta=b_1^2+b_2^2+b_3^2$. One easily finds that $\|B\|\le 1$ iff 
\begin{eqnarray}
b_2^2&\le&(1-b_1^2)(1-b_3^3)\cr
b_1^2&\le&1\cr
b_3^2&\le&1\label{condbi}
\end{eqnarray}
Since we can clearly suppose $b_1,b_2$ and $b_3$ to be nonnegative, and that $b_2$ reaches the maximal value allowed by \eqref{condbi}, we can set $b_1=\sin\alpha$, $b_3=\sin\beta$ and $b_2=\cos\alpha\cos\beta$ and look for the maximum of
\begin{equation}
\tilde f(\alpha,\beta)=w_1\sin\alpha+w_3\sin\beta+w_2\cos\alpha\cos\beta\label{tildef}
\end{equation}
for $\alpha,\beta\in[0,\pi/2]$. The gradient of $\tilde{f}$ vanishes when
\begin{eqnarray}
w_1\cos\alpha&=&w_2\sin\alpha\cos\beta\cr
w_3\cos\beta&=&w_2\cos\alpha\sin\beta\label{sysb}
\end{eqnarray}
Multiplying these two equations, we obtain, if $\cos\alpha\cos\beta\not=0$ :
\begin{equation}
w_1w_3=w_2^2\sin\alpha\sin\beta.\label{condw}
\end{equation}
Hence, $\tilde f$ has no critical point in the interior of $[0,\pi/2]^2$ when $w_1w_3\ge w_2^2$. In this case we inspect the values of $\tilde f$ on the boundary and see that its maximum is $w_1+w_3$. If $w_1w_3<w_2^2$ we can check that the critical point corresponds to a local minimum. Raising each equation of \eqref{sysb} and using \eqref{condw} readily yields $\cos^2\alpha$ and $\cos^2\beta$. We finally obtain at the critical point $\tilde f(\alpha,\beta)=\frac{\sqrt{w_1^2+w_2^2}\sqrt{w_3^2+w_2^2}}{w_2}$. We thus obtain the following result:
\begin{itemize}
\item If $w_2>\sqrt{w_1w_3}$, $\lambda(C)=\frac{\sqrt{w_1^2+w_2^2}\sqrt{w_3^2+w_2^2}}{w_2}$.
\item If $w_2\le \sqrt{w_1w_3}$, $\lambda(C)=w_1+w_3$.
\end{itemize}
We can check the continuity at $w_2=\sqrt{w_1w_3}$. We see that $w_2$ has no effect on the noncommutative length of $C$ as soon as it is smaller than the geometric mean of $w_1$ and $w_3$.

\bibliographystyle{unsrt}
\bibliography{../generalbib/SSTbiblio}
\end{document}